\documentclass[a4paper,12pt]{amsart}
\usepackage{amsfonts}
\usepackage{amssymb}
\usepackage{ifthen}
\usepackage{graphicx}

\newtheorem{thm}{Theorem}
\newtheorem{cor}[equation]{Corollary}
\newtheorem{lem}{Lemma}
\newtheorem{prop}[equation]{Proposition}

\newtheorem{conj}[equation]{Conjecture}
\newtheorem{rem}{Remark}
\theoremstyle{definition}
\newtheorem{defn}{Definition}

\newtheorem{prob}[equation]{Problem}
\newtheorem{ques}[equation]{Question}

\newcounter {own}
\def\theown {\thesection       .\arabic{own}}

\newenvironment{pf}[1][]{%
 \vskip 3mm
 \noindent
 \ifthenelse{\equal{#1}{}}%
  {{\slshape Proof. }}%
  {{\slshape #1.} }%
 }%
{\qed\bigskip}

\newcounter{alphabet}
\newcounter{tmp}

\makeatletter
\newcommand{\Ref}[1]{\@ifundefined{r@#1}{}{\setcounter{tmp}{\ref{#1}}\Alph{tmp}}}
\makeatother

\newenvironment{Lem}[1][]{\refstepcounter{alphabet}%
\bigskip%
\noindent%
{\bf Lemma \Alph{alphabet}}%
{\bf .} \itshape}{\vskip 8pt}

\newcommand{\ID}{{\mathbb D}}
\newcommand{\IB}{{\mathbb B}}

\def\be{\begin{equation}}
\def\ee{\end{equation}}

\newcommand{\bee}{\begin{enumerate}}
\newcommand{\eee}{\end{enumerate}}

\newcommand{\blem}{\begin{lem}}
\newcommand{\elem}{\end{lem}}
\newcommand{\bthm}{\begin{thm}}
\newcommand{\ethm}{\end{thm}}
\newcommand{\bcor}{\begin{cor}}
\newcommand{\ecor}{\end{cor}}
\newcommand{\beg}{\begin{examp}}
\newcommand{\eeg}{\end{examp}}
\newcommand{\begs}{\begin{examples}}
\newcommand{\eegs}{\end{examples}}
\newcommand{\bdefe}{\begin{defn}}
\newcommand{\edefe}{\end{defn}}
\newcommand{\bprob}{\begin{prob}}
\newcommand{\eprob}{\end{prob}}
\newcommand{\bques}{\begin{ques}}
\newcommand{\eques}{\end{ques}}
\newcommand{\bei}{\begin{itemize}}
\newcommand{\eei}{\end{itemize}}

\newcommand{\bcon}{\begin{conj}}
\newcommand{\econ}{\end{conj}}
\newcommand{\bcons}{\begin{conjs}}
\newcommand{\econs}{\end{conjs}}
\newcommand{\bprop}{\begin{prop}}
\newcommand{\eprop}{\end{prop}}
\newcommand{\br}{\begin{rem}}
\newcommand{\er}{\end{rem}}
\newcommand{\brs}{\begin{rems}}
\newcommand{\ers}{\end{rems}}
\newcommand{\bo}{\begin{obser}}
\newcommand{\eo}{\end{obser}}
\newcommand{\bos}{\begin{obsers}}
\newcommand{\eos}{\end{obsers}}
\newcommand{\bpf}{\begin{pf}}
\newcommand{\epf}{\end{pf}}
\newcommand{\ba}{\begin{array}}
\newcommand{\ea}{\end{array}}
\newcommand{\beq}{\begin{eqnarray}}
\newcommand{\beqq}{\begin{eqnarray*}}
\newcommand{\eeq}{\end{eqnarray}}
\newcommand{\eeqq}{\end{eqnarray*}}

\newcommand{\ra}{\rightarrow}

\newcommand{\ds}{\displaystyle}

\begin{document}
\bibliographystyle{amsplain}
\title [] {Equivalent moduli of continuity, Bloch's theorem for pluriharmonic mappings in $\mathbb{B}^{n}$}
\author{SH. Chen}
\address{Sh. Chen, Department of Mathematics,
Hunan Normal University, Changsha, Hunan 410081, People's Republic
of China.} \email{shlchen1982@yahoo.com.cn}

\author{S. Ponnusamy}
\address{S. Ponnusamy, Department of Mathematics,
Indian Institute of Technology Madras, Chennai-600 036, India.}
\email{samy@iitm.ac.in}
\author{X. Wang${}^{~\mathbf{*}}$}
\address{X. Wang, Department of Mathematics,
Hunan Normal University, Changsha, Hunan 410081, People's Republic
of China.} \email{xtwang@hunnu.edu.cn}

\subjclass[2000]{Primary: 30C65, 30C45; Secondary: 30C20}
\keywords{Pluriharmonic  mapping, Lipschitz-type space, Bloch constant,
Schwarz' Lemma, equivalent norm.\\
${}^{\mathbf{*}}$ Corresponding author}
\date{\today  
File: Ch-W-S12${}_{}$equiv-mod${}_{}$submit.tex}
\begin{abstract}
In this paper, we first establish a Schwarz-Pick type theorem for
pluriharmonic mappings and then we apply it to discuss the
equivalent norms on Lipschitz-type spaces. Finally, we obtain
several Landau's and  Bloch's type theorems for pluriharmonic
mappings.
\end{abstract}

\thanks{The research was partly supported by
NSF of China (No. 11071063). The work was carried out while the
first author was visiting IIT Madras, under ``RTFDCS Fellowship."
This author thanks Centre for International Co-operation in Science
(Formerly  Centre for Cooperation in Science \& Technology
among Developing Societies (CCSTDS)) for its support and cooperation.
}

\maketitle \pagestyle{myheadings} \markboth{SH. Chen,  S. Ponnusamy, X. Wang}{Equivalent moduli of
continuity of pluriharmonic mapping in $\mathbb{B}^{n}$}

\section{Introduction and Main Results}\label{csw-sec1}

Let $\mathbb{C}^{n}$ denote the complex normed (Euclidean) space of
dimension $n$. For $z=(z_{1},\ldots,z_{n})\in \mathbb{C}^{n}$,
the conjugate of $z$, denoted by $\overline{z}$, is defined by
$\overline{z}=(\overline{z}_{1},\ldots, \overline{z}_{n} ). $ For
$z$ and $w=(w_{1},\ldots,w_{n})\in\mathbb{C}^{n}
$, we write
$$\langle z,w\rangle :=z\cdot w= \sum_{k=1}^nz_k\overline{w}_k\quad \mbox{and}\quad
|z|:={\langle z,z\rangle}^{1/2}=(|z_{1}|^{2}+\cdots+|z_{n}|^{2})^{ 1/2}.
$$
For $a=(a_{1},\ldots,a_{n})\in \mathbb{C}^n$, we set $\IB^n(a, r)=\{z\in \mathbb{C}^{n}:\, |z-a|<r\}.
$
Also, we use $\IB^n$ to denote the unit ball $\IB^n(0, 1)$ and let
$\mathbb{D}=\{z\in \mathbb{C}:\; |z|<1\}$.

A continuous complex-valued function $f$ defined on a domain
$G\subset\mathbb{C}^{n}$ is said to be {\it pluriharmonic} if
for each fixed $z\in G$ and $\theta\in\partial\mathbb{B}^{n}$,
the function $f(z+\theta\zeta)$ is harmonic in $\{\zeta:\; |\zeta|< d_{G}(z)\}$, where $d_{G}(z)$ denotes the distance from $z$
to the boundary $\partial G$ of $G$. It follows from \cite[Theorem 4.4.9]{R} that
a real-valued function $u$ defined on $G$ is pluriharmonic if and only if $u$ is the real part of a holomorphic function
on $G$. Clearly, a mapping $f:\,\mathbb{B}^{n}\ra \mathbb{C}$ is pluriharmonic
if and only if $f$ has a representation $f=h+\overline{g}$, where
$g$ and $h$ are holomorphic mappings. For a pluriharmonic mapping
$f:\,\mathbb{B}^{n}\ra \mathbb{C}$, we introduce the notation
$$\nabla f=(f_{z_{1}},\ldots,f_{z_{n}})~\mbox{ and }~\nabla
\overline{f}=(f_{\overline{z}_{1}},\ldots,f_{\overline{z}_{n}}).
$$

For a proper domain $G$ of $\mathbb{C}^{n}$, let $\mathcal{H}_{k}(G)$ denote
the class of all pluriharmonic mappings $f=h+\overline{g}$ defined from $G$ into $\mathbb{C}$ such that for
any $\theta\in\partial\mathbb{B}^{n}$,
$$|\nabla \overline{f}(z)\cdot\theta|\leq k|\nabla
f(z)\cdot\overline{\theta}|
$$
for $z\in G$, where $k\in(0,1)$ is a constant, and both $h$ and $g$ are holomorphic in $G$.

Let $f$ be a sense-preserving  harmonic mapping from $\mathbb{D}$
into $\mathbb{C}$. We say that $f$ is a {\it $K$-quasiregular
harmonic mapping} if and only if
$$\frac{\Lambda_{f}(z)}{\lambda_{f}(z)}\leq K,  ~\mbox{ i.e., }~\frac{|f_{\overline{z}}(z)|}{|f_{z}(z)|}\leq \frac{K-1}{K+1}
$$
for $z\in\mathbb{D}$, where $\Lambda_{f}=|f_{z}|+|f_{\overline{z}}|$ and
$\lambda_{f}=|f_{z}|-|f_{\overline{z}}|.$

%


First we improve the Schwarz-Pick type theorem for $K$-quasiregular harmonic mappings obtained
recently by  Chen \cite[Theorem 7]{CH}.

\begin{lem}\label{lemma2.1}
Let $f$ be a sense-preserving and $K$-quasiregular harmonic mapping
on $\mathbb{D}$ with $f(\mathbb{D})\subset\mathbb{D}$. Then
\be\label{eq2.1}
\Lambda_{f}(z)\leq K\frac{1-|f(z)|^{2}}{1-|z|^{2}}\leq\frac{4K}{\pi}\left(\frac{\cos(|f(z)|\pi/2)}{1-|z|^{2}}\right), \quad
\mbox{$z\in\mathbb{D}$}.
\ee
Moreover, the first inequality of {\rm (\ref{eq2.1})} is sharp when $K=1.$
\end{lem}

By using Lemma \ref{lemma2.1}, we obtain a Schwarz-Pick type theorem
for pluriharmonic mappings which is as follows.

\begin{thm}\label{thm-schwarz-pick}
Let $f\in\mathcal{H}_{k}(\mathbb{B}^{n})$ and $|f(z)|<1$ for $z\in\mathbb{B}^{n}$. Then for each
$\theta\in\partial\mathbb{B}^{n}$,
$$|\nabla \overline{f}(z)\cdot\theta|+|\nabla f(z)\cdot\overline{\theta}|\leq
K \frac{1-|f(z)|^{2}}{1-|z|^{2}}, \quad K=\frac{1+k}{1-k}.
$$
\end{thm}

Proofs of Lemma \ref{lemma2.1} and Theorem \ref{thm-schwarz-pick}
will be given in Section \ref{csw-sec2}.

A continuous increasing function $\omega:\ [0,\infty)\rightarrow [0,\infty)$ with
$\omega(0)=0$ is called a {\it majorant} if $\omega(t)/t$ is
non-increasing for $t>0$. Given a subset $G$ of
$\mathbb{C}^{n}$, a function $f:\ G\rightarrow \mathbb{C}$ is
said to belong to the {\it Lipschitz space
$\Lambda_{\omega}(G)$} if there is a positive constant $C$ such
that
\be\label{eq1}
|f(z)-f(w)|\leq C\omega(|z-w|)
\ee
for all $z, w\in G.$ For $\delta_{0}>0$, let
\be\label{eq2}
\int_{0}^{\delta}\frac{\omega(t)}{t}\,dt\leq C\omega(\delta),\
0<\delta<\delta_{0}
\ee
and
\be\label{eq3}
\delta\int_{\delta}^{+\infty}\frac{\omega(t)}{t^{2}}\,dt\leq
C\omega(\delta),\ 0<\delta<\delta_{0}.
\ee

A majorant $\omega$ is said to be {\it regular} if it
satisfies the conditions (\ref{eq2}) and (\ref{eq3}) (see \cite{D,P}).






Let $G$ be a proper subdomain of $\mathbb{C}^{n}$. We say that a
function $f$ belongs to the {\it local Lipschitz space }
$\mbox{loc}\Lambda_{\omega}(G)$ if there is a constant $C>0$
satisfying (\ref{eq1}) for all $z$, $w\in G$ with
$|z-w|<\frac{1}{2}d_G(z)$. Moreover, $G$ is said to be a {\it
$\Lambda_{\omega}$-extension domain} if
$\Lambda_{\omega}(G)=\mbox{loc}\Lambda_{\omega}(G).$ The geometric
characterization of $\Lambda_{\omega}$-extension domains was first
given by Gehring and Martio \cite{GM}. Later, Lappalainen \cite{L}
extended it to the general case, and proved that $G$ is a
$\Lambda_{\omega}$-extension domain if and only if each pair of
points $z,w\in G$ can be joined by a rectifiable curve
$\gamma\subset G$ satisfying
\be\label{eq1.0}
\int_{\gamma}\frac{\omega(d_G(z))}{d_G(z)}\,ds(z) \leq C\omega(|z-w|)
\ee
with some fixed positive constant
$C=C(G,\omega)$, where $ds$ stands for the arclength measure on
$\gamma$. Furthermore, Lappalainen \cite[Theorem 4.12]{L} proved
that $\Lambda_{\omega}$-extension domains  exist only for majorants
$\omega$ satisfying the inequality (\ref{eq2}).

For $z_{1}$, $z_{2}\in G\subset\mathbb{C}^{n}$, let
$$d_{\omega,G}(z_{1},z_{2}):=\inf\int_{\gamma}\frac{\omega(d_G(z))}{d_G(z)}\,ds(z),
$$
where the infimum is taken over all rectifiable curves $\gamma\subset G$ joining $z_{1}$ and $z_{2}$. We say that
$f\in\Lambda_{\omega,\inf}(G)$ whenever
$$|f(z_{1})-f(z_{2})|\leq Cd_{\omega,G}(z_{1},z_{2}) ~\mbox{ for $z_{1},z_{2}\in G$},
$$
where $C$ is a positive constant which depends only on $f$ (see
\cite{KW}).

Dyakonov \cite{D} characterized the holomorphic functions in
$\Lambda_{\omega}$ in terms of their modulus. Later in
\cite[Theorems A and B]{P}, Pavlovi\'{c} came up with a relatively
simple proof of the results of Dyakonov. Recently, many authors
considered this topic and generalized the work of Dyakonov to
pseudo-holomorphic functions and real harmonic functions of several
variables for some special majorants $\omega(t)= t^{\alpha}$, where
$\alpha>0$ (see \cite{ABM,D1,KV,KM,KP,M1,M2,M3,MM}). By applying
Theorem \ref{thm-schwarz-pick}, we extend \cite[Theorems A and B]{P}
to the case of pluriharmonic mappings.

\begin{thm}\label{thm1}
Let $\omega$ be a  majorant satisfying {\rm (\ref{eq2})}, and let
$G$ be a $\Lambda_{\omega}$-extension. If $f\in\mathcal{H}_{k}(G)$
and is continuous up to the boundary $\partial G$, then
$$f\in\Lambda_{\omega}(G)\Longleftrightarrow |f|\in\Lambda_{\omega}(G)
\Longleftrightarrow |f|\in\Lambda_{\omega}(G,\partial G),
$$
where $\Lambda_{\omega}(G,\partial G)$ denotes the class of
continuous functions $f$ on $G\cup\partial G$ which satisfy {\rm
(\ref{eq1})} with some positive constant $C$, whenever $z\in G$ and
$w\in\partial G$.
\end{thm}

\begin{thm}\label{thm2}
Let $\omega$ be a  majorant satisfying {\rm (\ref{eq2})}. If
$f\in\mathcal{H}_{k}(G)$, then
$$f\in\Lambda_{\omega,\inf}(G)\Longleftrightarrow |f|\in\Lambda_{\omega,\inf}(G).
$$
\end{thm}

We remark that Theorems \ref{thm1} and \ref{thm2} are  the generalizations of \cite[Theorem 1]{CPW} and
\cite[Theorem 2]{CPW}, respectively.

To state our final result, we need some preparations. First we
recall that a mapping $f:\,\Omega \ra \mathbb{C}^{n}$ is said
to be {\it vector-valued pluriharmonic} if every component of $f$ is
pluriharmonic. Let $H(\IB^n, \mathbb{C}^n)$ denote the set of all
pluriharmonic mappings from $\mathbb{B}^{n}$ into $\mathbb{C}^{n}$.
Obviously, a mapping $f\in H(\IB^n, \mathbb{C}^n)$  is pluriharmonic
if and only if $f$ has a representation $f=h+\overline{g}$, where
$g$ and $h$ are holomorphic mappings $\mathbb{B}^{n}$ into
$\mathbb{C}^{n}$. It is convenient to identify each point
$z=(z_{1},\ldots,z_{n})\in\mathbb{C}^{n}$ with an $n\times 1$ column
matrix so that
$$z=\left(\begin{array}{cccc}
z_{1}   \\
\vdots \\
 z_{n}
\end{array}\right), \quad dz=\left(\begin{array}{cccc}
dz_{1}   \\
\vdots \\
 dz_{n}
\end{array}\right)~ \mbox{ and }~d\overline{z}=\left(\begin{array}{cccc}
d\overline{z}_{1}   \\
\vdots \\
 d\overline{z}_{n}
\end{array}\right).$$
For a  $f=(f_{1},\ldots,f_{n})\in H(\IB^n, \mathbb{C}^n)$, we denote
by $\partial f/\partial z_{j}$ the column vector formed by $\partial
f_{1}/\partial z_{j},\ldots ,\partial f_{n}/\partial z_{j}$, and
$$f_{z}=\left (\frac{\partial f}{\partial z_{1}}~\cdots ~\frac{\partial f}{\partial z_{n}}\right ) :=
\left ( \frac{\partial f_i}{\partial z_{j}}\right )_{n\times n},
$$
the $n\times n$ matrix formed by these column vectors, namely, by the complex gradients
$\nabla f_1,\ldots,\nabla f_n$.
Similarly,
$$f_{\overline{z}}=\left (\frac{\partial f}{\partial \overline{z}_{1}}~\cdots ~\frac{\partial f}{\partial
\overline{z}_{n}}\right )
 :=\left ( \frac{\partial f_i}{\partial \overline{z}_{j}}\right )_{n\times n}
$$
the $n\times n$ matrix formed by the column vectors $\partial f/\partial \overline{z}_{j}$ for $j\in\{1,\ldots , n\}$.
For an $n\times n$ matrix $A$, we introduce the operator norm
$$|A|=\sup_{x\neq 0}\frac{|Ax|}{|x|}=\max\{|A\theta|:\
\theta\in\partial \mathbb{B}^{n}\}.
$$
For pluriharmonic mappings $f:\,\mathbb{B}^{n}\ra \mathbb{C}^{n}$, we use the following standard notations (cf.
\cite{HG1}):
$$\Lambda_{f}(z)=\max_{
\theta\in\partial\mathbb{B}^{n}}|f_{z}(z)\theta+f_{\overline{z}}(z)\overline{\theta}|
\ \mbox{and}\ \lambda_{f}(z)=\min_{
\theta\in\partial\mathbb{B}^{n}}|f_{z}(z)\theta+f_{\overline{z}}(z)\overline{\theta}|.
$$

Let $f=(f_{1},\ldots,f_{n})\in H(\IB^n, \mathbb{C}^n)$.
For $j\in\{1,\ldots,n\}$, we let
$z=(z_{1},\ldots,z_{n})$, $z_{j}=x_{j}+iy_{j}$ and
$f_{j}(z)=u_{j}(z)+iv_{j}(z)$, where $u_{j}$ and $v_{j}$ are real
pluriharmonic functions from $\mathbb{B}^{n}$ into $\mathbb{R}$. We
denote the real Jacobian matrix of $f$ by
$$J_{f}=\left(\begin{array}{cccc}
\ds \frac{\partial u_{1}}{\partial x_{1}}\; \frac{\partial
u_{1}}{\partial y_{1}}\; \frac{\partial u_{1}}{\partial x_{2}}\;
\frac{\partial u_{1}}{\partial y_{2}}\;\cdots\;
 \frac{\partial u_{1}}{\partial x_{n}}\; \frac{\partial u_{1}}{\partial y_{n}}\\[4mm]
 \ds \frac{\partial v_{1}}{\partial x_{1}}\; \frac{\partial v_{1}}{\partial y_{1}}\;
\frac{\partial v_{1}}{\partial x_{2}}\; \frac{\partial
v_{1}}{\partial y_{2}}\;\cdots\;
 \frac{\partial v_{1}}{\partial x_{n}}\; \frac{\partial v_{1}}{\partial
 y_{n}}\\[2mm]
\vdots \\[2mm]
 \ds \frac{\partial u_{n}}{\partial x_{1}}\; \frac{\partial u_{n}}{\partial y_{1}}\;
\frac{\partial u_{n}}{\partial x_{2}}\; \frac{\partial
u_{n}}{\partial y_{2}}\;\cdots\;
 \frac{\partial u_{n}}{\partial x_{n}}\; \frac{\partial u_{n}}{\partial y_{n}}\\[4mm]
\ds  \frac{\partial v_{n}}{\partial x_{1}}\; \frac{\partial v_{n}}{\partial y_{1}}\;
\frac{\partial v_{n}}{\partial x_{2}}\; \frac{\partial
v_{n}}{\partial y_{2}}\;\cdots\;
 \frac{\partial v_{n}}{\partial x_{n}}\; \frac{\partial v_{n}}{\partial
 y_{n}}
\end{array}\right).
$$
Let $\mathbb{B}^{2n}_{\mathbb{R}}$ denote the unit ball of
$\mathbb{R}^{2n}$. Then (see \cite{HG1})
$$\Lambda_{f}=\max_{\theta\in\partial\mathbb{B}^{2n}_{\mathbb{R}}}|J_{f}\theta|\
\mbox{and}\
\lambda_{f}=\min_{\theta\in\partial\mathbb{B}^{2n}_{\mathbb{R}}}|J_{f}\theta|
$$

We use $b^{p}_{h}(\mathbb{B}^{n},\mathbb{C}^{n})$ to denote the {\it
pluriharmonic Bergman space} consisting of all pluriharmonic
mappings $f\in H(\IB^n, \mathbb{C}^n)$ such that
$$\|f\|_{b^{p}}=\left (\int_{\mathbb{B}^{n}}|f(z)|^{p}
\,dV(z)\right )^{1/p} <\infty\ \mbox{or}\ \|f\|_{b_{N}^{p}}=\left
(\int_{\mathbb{B}^{n}}|f(z)|^{p} \,dV_{N}(z)\right )^{1/p} <\infty,
$$
where $p\in(0,\infty)$, $n\geq2,$ $dV$ denotes the Lebesgue volume
measure on $\mathbb{C}^{n}$ and $dV_{N}$ denotes the normalized
Lebesgue volume measure on $\mathbb{B}^{n}$. Obviously, if $f\in
H(\IB^n, \mathbb{C}^n) $ and $f$ is bounded, then $f\in
b^{p}_{h}(\mathbb{B}^{n},\mathbb{C}^{n}).$

\begin{thm}\label{thm8}
 Let $r\in(0,1)$ and $f\in b^{p}_{h}(\mathbb{B}^{n},\mathbb{C}^{n})$ with $\|f\|_{b_{N}^{p}}\leq M$, $f(0)=0$ and $\det J_{f}(0)=\alpha>0$.
 Then $f$ is injective in $\mathbb{B}^{n}(0,r\rho(r))$ with
$$\rho(r)=\frac{\alpha\pi^{2n+1}}{4m(4M(r))^{2n}}$$
and $f(\mathbb{B}^{n}(0,r\rho(r)))$ contains a univalent ball with
the radius
$$R\geq\max_{0<r<1}\left\{\frac{\alpha\pi^{4n}r}{8m(4M(r))^{4n-1}}\right\},$$
where
\be\label{extraeq2}
M(r)=\frac{M}{r(1-r)^{2n/p}}~\mbox{ and }~m=2\sqrt{2}\left
( \frac{3+\sqrt{17}}{(\sqrt{5-\sqrt{17}})(1+\sqrt{17})}\right
)\approx 4.2.
\ee
\end{thm}

We remark that Theorem \ref{thm8} is a generalization of \cite[Theorem 5]{HG1}.
We now recall that a holomorphic function $f:\ \mathbb{B}^{n}\rightarrow\mathbb{C}^{n}$
is {\it convex} in $\mathbb{B}^{n}$ if it is one-to-one and the range $f(\mathbb{B}^{n})$ is a
convex domain.

\begin{thm}\label{thm5}
Suppose $f=h+\overline{g}\in H(\IB^n, \mathbb{C}^n)$, $f(0)=0,$ $| f_{\overline{z}}(0)|=0$
and $\det f_{z}(0)=I_{n}$, where $h$ is a convex biholomorphic
mapping and $g$ is a holomorphic mapping. If for any
$z\in\mathbb{B}^{n}$, $|f_{\overline{z}}(z)|\leq |f_{z}(z)|$, then
$f$ is univalent in $\mathbb{B}^{n}(0,\rho_{1})$, where
$$\rho_{1}=\frac{1}{m_{2}+m_{3}} ~\mbox{ with $m_2\approx 9.444$ and $m_3=6.75$.}
$$
Moreover, the range $f(\mathbb{B}^{n}(0,\rho_{1}))$ contains a univalent ball with
center $0$ and radius at least $R_{1}$, where
$R_{1}=\frac{\rho_{1}}{2}.$
\end{thm}

The precise values of $m_2$ and $m_3$ are given in the proof of Theorem \ref{thm5}.

A continuous mapping $f:\ \Omega\subset\mathbb{R}^{n}\rightarrow\mathbb{R}^{n}$ is called {\it quasiregular} if $f\in
W_{n,{\rm loc}}^{1}(\Omega)$ and
$$|f'(x)|^{n}\leq KJ_{f}(x)\ \mbox{for almost every}\ x\in\Omega,
$$
where $K$ $(\geq1)$ is a constant, $f\in W_{n,{\rm loc}}^{1}(\Omega)$ means that the
distributional derivatives $\partial f_{j}/\partial x_{k}$ of the
coordinates $f_{j}$ of $f $ are locally in $L^{n}$ and $J_{f}(x)$
denotes the Jacobian of $f$ (cf. \cite{V}).

\begin{defn}\label{def1.1}
A pluriharmonic mapping $f:\,\mathbb{B}^{n}\ra \mathbb{C}^{n}$ is said to be a {\it $(K,K_{1})$-pluriharmonic
mapping} if for each $z\in\mathbb{B}^{n}$ and $\theta\in\partial\mathbb{B}^{n}$,
\be\label{eqth3}
|f_{z}(z)|^{n}\leq K|\det f_{z}(z)|\ \mbox{ and }\ K_{1}|
f_{\overline{z}}(z)\theta|\leq | f_{z}(z)\overline{\theta}|,
\ee
where $K$ $(\geq 1)$ and $ K_{1}$ $(>1)$ are constants.
\end{defn}

Obviously, every ($K,K_{1}$)-pluriharmonic mapping $f:\,\mathbb{B}^{n}\ra \mathbb{C}^{n}$ is
called Wu $K$-mapping if $f_{\overline{z}}\equiv 0$ (see \cite{W}). In fact, holomorphic $K$-quasiregular mappings are
referred to as Wu $K^{1-\frac{1}{n}}$-mappings (cf. \cite{HG,W}).

For a holomorphic mapping $f$ from the unit ball $\mathbb{B}^{n}$
into $\mathbb{C}^{n}$, $\mathbb{B}^{n}(a,r)$ is called a {\it
schlicht ball} of $f$ if there is a subregion
$\Omega\subset\mathbb{B}^{n}$ such that $f$ maps $\Omega$
biholomorphically onto $\mathbb{B}^{n}(a,r)$. We denote by $B_{f}$
the least upper bound of radii of all schlicht balls contained in
$f(\mathbb{B}^{n})$ and call this the Bloch radius of $f$. The
classical theorem of Bloch for holomorphic functions in the unit
disk fails to extend to general holomorphic mappings in the ball of
$\mathbb{C}^{n}$ (see \cite{T,W}). However, in 1946, Bochner
\cite{B} proved that Bloch's theorem does hold for a class of real
harmonic quasiregular mappings.  Recently, Chen and Gauthier
\cite{HG1} proved that  Bloch's theorem also holds for a class of
pluriharmonic $K$-mappings.


In this paper, our last aim is to prove the existence of Bloch's
constant for a new class of pluriharmonic mappings. Our result is
also a generalization of \cite[Theorem 6]{HG}.
 We now state a version of
Bloch's theorem for a class of ($K,K_{1}$)-quasiregular
pluriharmonic mappings.

\begin{thm}\label{thm3}
Suppose $f$ is a $(K,K_{1})$-quasiregular pluriharmonic mapping of
$\mathbb{B}^{n}$ into $\mathbb{C}^{n}$ with  $|J_{f}(0)|=1$. Then
$f(\mathbb{B}^{n})$ contains a schlicht ball with radius at least
$$B_{f}\geq\max_{0<t<1}\left\{\frac{\pi^{4n}t}{8m(4M(t))^{4n-1}}\right\},\quad
\mbox{with }~
M(t)=\frac{K^{\frac{1}{n}}(1+K_{1})}{tK_{1}}\log\left(\frac{1}{1-t}\right),
$$
where $m$  defined as in Theorem {\rm \ref{thm8}}.
\end{thm}

Proofs of Theorems \ref{thm1} and \ref{thm2} will be given in
Section \ref{csw-sec3} while the proofs of Theorems \ref{thm8},
\ref{thm5} and \ref{thm3}  in Section \ref{csw-sec4}.

\section{ Schwarz-Pick lemma for pluriharmonic mappings in $\mathbb{B}^{n}$}\label{csw-sec2}

Let   $\Omega$   be a domain in $\mathbb{C}$ and $\rho>0$
 a conformal metric in   $\Omega$. The Gaussian curvature
of the domain is given by $K_{\rho}=-(1/(2\rho))\Delta\log\rho.$
We denote by $\lambda(z)|dz|^{2}$ the hyperbolic metric in $\mathbb{D}$, where
$\lambda(z)=4/(1-|z|^{2})^{2}$.


\begin{Lem}{\rm \bf (Ahlfors-Schwarz lemma)}\label{AL}
If $\rho>0$ is a $C^{2}$-function $($metric density$)$ in
$\mathbb{D}$ and Gaussian curvature $K_{\rho}\leq-1$, then
$\rho\leq\lambda$ $($cf. \cite{A}$)$.
\end{Lem}




\subsection*{Proof of Lemma \ref{lemma2.1}}
By assumption, we observe that $f$ is an open mapping, and so $|f_{z}(z)|\neq 0$
in $\ID$. Let
$$\rho(z)=\frac{4}{(K+1)^{2}}\lambda(f(z))|f_{z}(z)|^{2}, \quad z\in \ID.
$$
Then $ds^{2}=\rho(z)|dz|^{2}$.  Simple calculations yield
\begin{eqnarray*}
\Delta\log\rho(z)&=&\Delta\log\left [\frac{4}{(K+1)^{2}}\lambda(f(z))|f_{z}(z)|^{2}\right ]\\
&=&4\left (\log(\lambda(f(z)))\right )_{z\overline{z}}\\
&=&\frac{8|f_{z}(z)|^{2}}{(1-|f(z)|^{2})^{2}}\left [1+\frac{|f_{\overline{z}}(z)|^{2}}{|f_{z}(z)|^{2}}
+2\mbox{Re}\left (\frac{f^{2}(z)\overline{f_{z}}(z)\overline{f_{\overline{z}}}(z)}{|f_{z}(z)|^{2}}\right )\right ]\\
&=&\frac{(K+1)^{2}\rho(z)}{2}\left [1+\frac{|f_{\overline{z}}(z)|^{2}}{|f_{z}(z)|^{2}}
+2\mbox{Re}\left(\frac{f^{2}(z)\overline{f_{z}}(z)\overline{f_{\overline{z}}}(z)}{|f_{z}(z)|^{2}}\right)\right ]\\
&\geq & \frac{(K+1)^{2}\rho(z)}{2}\left  ( 1-\frac{|f_{\overline{z}}(z)|}{|f_{z}(z)|}\right )^2\\
&\geq & \frac{(K+1)^{2}\rho(z)}{2}\left  ( 1-\frac{K-1}{K+1}\right )^2
=2\rho(z)
\end{eqnarray*}
which, together with the definition of $K_{\rho}$, gives $ K_{\rho}(z) \leq -1$.
Thus, by Lemma \Ref{AL}, we have
$$\rho(z)=\frac{4}{(K+1)^{2}}\lambda(f(z))|f_{z}(z)|^{2}\leq\lambda(z)
$$
whence
$$\Lambda_{f}(z)\leq\frac{2K}{1+K}|f_{z}(z)|\leq K\frac{1-|f(z)|^{2}}{1-|z|^{2}}.
$$
The proof of the lemma is complete. \qed

\subsection*{Proof of Theorem \ref{thm-schwarz-pick}}
For each fixed $\theta\in\partial\mathbb{B}^{n}$, let $F(\zeta)=f(\theta\zeta)$ in $\mathbb{D}.$
Then $F$ is harmonic and $|F(\zeta)|< 1$ on $\mathbb{D}$. It follows that
$$\Lambda _F=|F_{\zeta}|+|F_{\overline{\zeta}}|=|\nabla f\cdot\overline{\theta}|
+|\nabla \overline{f}\cdot\theta| \leq K(|\nabla f\cdot\overline{\theta}|-|\nabla
\overline{f}\cdot\theta|)=K(|F_{\zeta}|-|F_{\overline{\zeta}}|)
$$
which implies that $F$ is a $K$-quasiregular harmonic mapping in
$\mathbb{D}$, where $ K=\frac{1+k}{1-k}.$ Hence, Lemma \ref{lemma2.1} shows that
$$ |\nabla f(z)\cdot\overline{\theta}|+|\nabla \overline{f}(z)\cdot\theta| =
\Lambda _F(\zeta) \leq
K\frac{1-|F(\zeta)|^{2}}{1-|\zeta|^{2}}=K\frac{1-|f(z)|^{2}}{1-|z|^{2}},
$$
where $z=\zeta\theta.$ This completes the proof. \qed

\section{Equivalent moduli of continuity  for pluriharmonic mappings }\label{csw-sec3}

\subsection*{Proof of Theorem \ref{thm1}}
The implications
$``f\in\Lambda_{\omega}(G)\Rightarrow |f|\in\Lambda_{\omega}(G)\Rightarrow
|f|\in\Lambda_{\omega}(G,\partial G)"
$
are obvious. Therefore, we only need to prove the implication:
$|f|\in\Lambda_{\omega}(G,\partial G)\Rightarrow
f\in\Lambda_{\omega}(G).$ In order to prove this, for a fixed $z\in G$, we let
\be\label{extraeq1}
M_{z}:=\sup\{|f(\zeta)|:\ |\zeta-z|<d_G(z)\},
\ee
and define the following function:
$$F(\eta)=\frac{f(z+d_G(z)\eta)}{M_{z}}, \quad  \eta\in \mathbb{B}^{n}.
$$
By a simple calculation, we obtain that for $\theta\in\partial\mathbb{B}^{n}$,
$$\left |\nabla \overline{F}(\eta)\cdot\theta\right |
=\frac{d_G(z)}{M_{z}} \left |\nabla \overline{f}(\xi)\cdot\theta\right
| \leq\frac{kd_G(z)}{M_{z}}\left |\nabla
f(\xi)\cdot\overline{\theta}\right | =k \left |\nabla
F(\eta)\cdot\overline{\theta}\right |
$$
where $\xi=z+d_G(z)\eta.$ Then, $F\in\mathcal{H}_{k}(\mathbb{B}^{n})$
and $|F(\eta)|\leq1$ in $\mathbb{B}^{n}$. By Theorem \ref{thm-schwarz-pick}, we have that for
$\theta\in\partial\mathbb{B}^{n}$,
$$|\nabla F(0)\cdot\overline{\theta}|+|\nabla \overline{F}(0)\cdot\theta|\leq K(1-|F(0)|^{2})
$$
which in turn gives
\be\label{eq4}
d_G(z)\left (|\nabla f(z)\cdot\overline{\theta}|+|\nabla
\overline{f}(z)\cdot\theta|\right) \leq 2K(M_{z}-|f(z)|), \quad K=\frac{1+k}{1-k}.
\ee
For a fixed $\varepsilon_{0}>0$, there exists a $\zeta\in\partial G$
such that $|\zeta-z|<(1+\varepsilon_{0})d_G(z)$.  Then, for
$w\in\mathbb{B}^{n}(z,d_G(z))$, we have
\begin{eqnarray*}
|f(w)|-|f(z)|&\leq&\big ||f(w)|-|f(\zeta)|\big |+\big ||f(\zeta)|-|f(z)|\big |\\
&\leq&C\omega((2+\varepsilon_{0})d_G(z))+C\omega((1+\varepsilon_{0})d_G(z)),
\end{eqnarray*}
where $C$ is a positive constant. Now we take $\varepsilon_{0}=1$.
Then
$$\sup_{w\in\mathbb{B}^{n}(z,d_G(z))}(|f(w)|-|f(z)|)\leq |f(w)|-|f(z)|\leq5C\omega(d_G(z))
$$
whence $M_{z}-|f(z)|\leq 5C\omega(d_G(z)),$
where $C$ is a positive constant. Thus for any $\theta\in\mathbb{B}^{n}$, by (\ref{eq4}) and the last inequality, we have
\be\label{eq4.00} |\nabla
f(z)\cdot\overline{\theta}|+|\nabla
\overline{f}(z)\cdot\theta|\leq10CK\cdot\frac{\omega(d_G(z))}{d_G(z)} ~\mbox{ for $z\in G$.}
\ee
For points $z_{1}, z_{2}\in G$, let
$\gamma\subset G$ be a rectifiable curve which joins $z_{1}$ and $z_{2}$
satisfying (\ref{eq1.0}). Integrating (\ref{eq4.00}) along $\gamma$,
we obtain that
\be\label{eq6}
|f(z_{1})-f(z_{2})|
\leq10CK\int_{\gamma}\frac{\omega(d_G(z))}{d_G(z)}\,ds(z).
\ee
Therefore, (\ref{eq1.0}) and (\ref{eq6}) yield
$|f(z_{1})-f(z_{2})|\leq C_{1}\cdot\omega(|z_{1}-z_{2}|),
$
where $C_{1}$ is a positive constant. This completes the proof.
\qed
\subsection*{Proof of Theorem \ref{thm2}} The implication $f\in\Lambda_{\omega,\inf}(G)\Rightarrow
|f|\in\Lambda_{\omega,\inf}(G)$ is obvious. We need only to prove
that $|f|\in\Lambda_{\omega,\inf}(G)\Rightarrow
f\in\Lambda_{\omega,\inf}(G)$.

Assume that $|f|\in\Lambda_{\omega,\inf}(G)$ and fix $z\in G$. Then it follows from a similar reasoning as
in the proof of the inequality (\ref{eq4}) that for $\theta\in\partial\mathbb{B}^{n}$,
\be\label{eq4.0}
d_G(z)(|\nabla f(z)\cdot\overline{\theta}|+|\nabla
\overline{f}(z)\cdot\theta|)\leq2K(M_{z}-|f(z)|),
\ee
where $M_z$ is defined by \eqref{extraeq1}.
For $w\in\mathbb{B}^{n}(z,d_G(z))$, there exists a positive constant $C$ such that
\be\label{eq10} |f(w)|-|f(z)|\leq C d_{\omega,G}(w,z)\leq
C\int_{[w,z]}\frac{\omega(d_G(\zeta))}{d_G(\zeta)}\,ds(\zeta),
\ee
where $[w,z]$ denotes the straight segment with endpoints $w$ and
$z$. We observe that if $\zeta\in [w,z]$, then one has
$[w,z]\subset\mathbb{B}^{n}(z,d_G(z))\subset G$ and therefore,
$$d_G(\zeta)\geq d_{\mathbb{B}^{n}(z,d_G(z))}(\zeta)
$$
which gives
\be\label{eq11}
\frac{\omega(d_G(\zeta))}{d_G(\zeta)}\leq\frac{\omega(d_{\mathbb{B}^{n}(z,d_G(z))}(\zeta))}
{d_{\mathbb{B}^{n}(z,d_G(z))}(\zeta)}.
\ee
For each $w\in\mathbb{B}^{n}(z,d_G(z))$, (\ref{eq10}) and (\ref{eq11}) imply that
\begin{eqnarray*}
|f(w)|-|f(z)|&\leq&C\int_{[w,z]}\frac{\omega(d_G(\zeta))}{d_G(\zeta)}\,ds(\zeta)\\
&\leq&C\int_{[w,z]}\frac{\omega(d_{\mathbb{B}^{n}(z,d_G(z))}(\zeta))}
{d_{\mathbb{B}^{n}(z,d_G(z))}(\zeta)}\,ds(\zeta)\\
&=&C\int_{[w,z]}\frac{\omega(d_G(z)-|\zeta-z|)}{d_G(z)-|\zeta-z|}\,ds(\zeta)\\
&\leq&C\int_{0}^{d_G(z)}\frac{\omega(t)}{t}\,dt\\
&\leq&C\omega(d_G(z)).
\end{eqnarray*}
From the last inequality,  we obtain that
\be\label{eq14}
M_{z}-|f(z)|\leq C\omega(d_G(z)).
\ee
Again, for any $z_{1}$, $z_{2}\in G$, by
(\ref{eq4.0}) and (\ref{eq14}),  there exists a positive constant
$C_{1}$ such that
$|f(z_{1})-f(z_{2})|\leq C_{1}d_{\omega,G}(z_{1},z_{2}).
$
The proof of the theorem is complete.
\qed

\section{Landau's and Bloch's theorem for pluriharmonic mappings }\label{csw-sec4}

The following three Lemmas are useful for the proof of Theorem
\ref{thm3}.

\begin{Lem}{\rm (\cite[Lemma 1]{HG1} or \cite[Lemma 4]{LX})}\label{LemB}
Let $A$ be an $n \times n$ complex $($real$)$ matrix. Then for any unit
vector $\theta\in\partial \mathbb{B}^{n}$, the inequality
$|A\theta|\geq|\det A|/|A|^{n-1}
$
holds.
\end{Lem}


\begin{Lem}{\rm (\cite[Lemma 4]{HG})}\label{LemD}
Let $A$ be a holomorphic mapping from $\mathbb{B}^{n}(0,r)$ into the
space of $n \times n$ complex matrices. If $A(0)=0$ and $|A(z)|\leq
M$ in $\mathbb{B}^{n}(0,r),$ then
$$|A(z)|\leq \frac{M|z|}{r}, \quad z\in \mathbb{B}^{n}(0,r).
$$
\end{Lem}

\subsection*{Proof of Theorem \ref{thm8}}
Fix $z\in\mathbb{B}^{n}$ and let $D_{z}=\{\zeta\in\mathbb{C}^{n}:\ |\zeta-z|<1-|z|\}.$
Then by Jensen's inequality, for $r\in[0,1-|z|)$ and $p\in[1,\infty)$, we have
\be\label{eqb5}
|f(z)|^{p}\leq\int_{\partial\mathbb{B}^{n}}|f(z+r\zeta)|^{p}\,d\sigma(\zeta).
\ee
Multiplying the formula (\ref{eqb5}) by $2nr^{2n-1}$ and integrating from $0$ to $1-|z|$, we have
\begin{eqnarray*}
(1-|z|)^{2n}|f(z)|^{p}&\leq&\int_{0}^{1-|z|}\left[2nr^{2n-1}\int_{\partial\mathbb{B}^{n}}
|f(z+r\zeta)|^{p}\,d\sigma(\zeta)\right]dr\\
&=&\int_{D_{z}}|f(z)|^{p}\,dV(z)\\
&\leq&\int_{\mathbb{B}^{n}}|f(z)|^{p}\,dV_{N}(z)
\leq M^{p}
\end{eqnarray*}
which gives
$$|f(z)|\leq\frac{M}{(1-|z|)^{2n/p}}.
$$
For $\zeta\in\mathbb{B}^{n}$ and $r\in(0,1)$, let $F(\zeta)=r^{-1}f(r\zeta)$. Then
$$|F(\zeta)|\leq\frac{M}{r(1-r)^{2n/p}}=M(r)\ \mbox{and}\ J_{F}(0)=J_{f}(0)=\alpha.
$$
Using \cite[Theorem 5]{HG1}, we obtain that $f$ is injective in
$\mathbb{B}^{n}(0,r\rho(r))$ with
$$\rho(r)=\frac{\alpha\pi^{2n+1}}{4m(4M(r))^{2n}}$$
and $f(\mathbb{B}^{n}(0,r\rho(r)))$ contains a univalent ball with radius
$$R\geq\max_{0<r<1}\left\{\frac{\alpha\pi^{4n}r}{8m(4M(r))^{4n-1}}\right\},
$$
where $m$ is given by \eqref{extraeq2}.
The proof is complete. \qed

For the proof of Theorem \ref{thm5}, we need the following lemma due to Fitzgerald and Thomas
is from \cite{FT}.

\begin{Lem}$($\cite[Proposition 2.2]{FT}$)$\label{ThmH}
Let $f$ be a convex mapping from $\mathbb{B}^{n}$ into
$\mathbb{C}^{n}$ with $f(0)=0$  and $f'(0)=I_{n}$, the $n\times n$
identity matrix. Suppose $t$ is a positive integer, $\theta \in
\partial \IB^n$ and $r\in (0,1)$. Then
$$|D_{\theta}^{t}f(r\theta)|\leq\frac{t!}{(1-r)^{t+1}}.$$
\end{Lem}

\subsection*{Proof of Theorem \ref{thm5}}
We begin to note that Lemma  \Ref{ThmH} gives
$$|f_{z}(z)-f_{z}(0)|\leq 1+\frac{1}{(1-|z|)^{2}} ~\mbox{ for $z\in\mathbb{B}^{n}$.}
$$
Let $W_{2}(r)=[1+(1-r)^{2}]/[r(1-r)^{2}]$ for $r\in (0,1)$. Then
$$ W_{2}(r_2)=\min_{r\in (0,1)}\{W_{2}(r)\},\quad \mbox{with }~ r_{2}=1-\sqrt[3]{1+\sqrt{2}}+\frac{1}{\sqrt[3]{1+\sqrt{2}}}\approx 0.404.
$$
Denote $W_{2}(r_2)$ by $m_2$. Then $m_2\approx 9.444$ and, by Lemma \Ref{LemD}, we have
$$|f_{z}(z)-f_{z}(0)|\leq m_{2}|z| ~\mbox{ for $|z|\leq r_{2}$,}
$$
By hypotheses, we have
$$|f_{\overline{z}}(z)-f_{\overline{z}}(0)|\leq |f_{z}(z)|\leq \frac{1}{
(1-|z|)^{2}} ~\mbox{  for $z\in\mathbb{B}^{n}$.}
$$
Let $W_{3}(r)=1/[r(1-r)^{2}]$ for $r\in (0,1)$. Then
$$W_{3}(r_3)=\min_{r\in (0,1)}\{W_{3}(r)\}, ~\mbox{ with }~ r_{3}=1/3.
$$
We denote $W_{3}(r_3)$ by $m_3$. Then $m_3=6.75$ and,  by Lemma \Ref{LemD}, we have
$$|f_{z}(z)-f_{z}(0)|\leq m_{3}|z|~\mbox{ for $|z|\leq r_3$.}
$$
Hence for $z\in\mathbb{B}^{n}(0,\rho_{1})$ with $\rho_{1}\leq r_{3}$,
$$|f_{z}(z)-f_{z}(0)|\leq m_{2}|z|\ \mbox{ and } \ |f_{\overline{z}}(z)-f_{\overline{z}}(0)|\leq m_{3}|z|.
$$

In order to prove the univalence of $f$ in
$\mathbb{B}^{n}(0,\rho_{1})$, we choose two distinct points $z'$,
$z''\in \mathbb{B}^{n}(0,\rho_{1})$
  and let $[z',z'']$ denote the segment
from $z'$ to $z''$ with the endpoints $z'$  and $z''$. Then
\begin{eqnarray*}
|f(z')-f(z'')|&\geq&
\left|\int_{[z',z'']}f_{z}(0)dz+f_{\overline{z}}(0)\,d\overline{z}\right|\\
&& -\left|\int_{[z',z'']}(f_{z}(z)-f_{z}(0))\,dz+(f_{\overline{z}}(z)-f_{\overline{z}}(0))\,d\overline{z}\right|\\
&\geq&|z'-z''|\big\{1-(m_{2}+m_{3})\rho_{1}\big\}\\
&>&0
\end{eqnarray*}
which shows that $f$ is univalent in $\mathbb{B}^{n}(0,\rho_{1})$. Furthermore, for any $z$ with $|z|=\rho_{1}$, we have
\begin{eqnarray*}
|f(z)-f(0)|&\geq&
\left|\int_{[0,\rho_{1}]}f_{z}(0)\,dz+f_{\overline{z}}(0)\,d\overline{z}\right|\\
&& -\left|\int_{[0,\rho_{1}]}(f_{z}(z)-f_{z}(0))\,dz+(f_{\overline{z}}(z)-f_{\overline{z}}(0))\,d\overline{z}\right|\\
&\geq&\frac{\rho_{1}}{2}.
\end{eqnarray*}
The proof of this theorem is complete. \qed

\subsection*{Proof of Theorem \ref{thm3}}
 Without loss of generality, we assume that $f$ is
pluriharmonic on $\overline{\mathbb{B}}^{n}$. Otherwise we replace
$f(z)$ by $f(sz)$ for some $s\in (0, 1)$. Then there exists some
$z_{0}\in\mathbb{B}^{n}$  such that
\begin{enumerate}
\item $(1-|z_{0}|)^{n}|\det f_{z}(z_{0})|=1$; and
\item $(1-|z|)^{n}|\det f_{z}(z)|\leq 1$ for all $z$ in the set $\{z:\,|z_{0}|=r\leq |z|\leq 1\}$.
\end{enumerate}
Hence it follows from the fact $|\det f_{z}(z)|\leq|\det f_{z}(z_{0})|$ for any $z$ in
$\{z:\; |z|=r=|z_0|\}$ and the maximum principle that
$$|\det f_{z}(z)|\leq|\det f_{z}(z_{0})|
$$
in the disk $\{z:\; |z|\leq r\}$. For $\zeta\in\mathbb{B}^{n}$ and $t\in(0,1)$, let
$$ F(\zeta)=\frac{1}{t}[f(p(\zeta))-f(z_{0})],
$$
where $p(\zeta)=z_{0}+t(1-r)\zeta$. Then
\be\label{eqn}
|\det F_{\zeta}(\zeta)|\leq\frac{1}{(1-t|\zeta|)^{n}}
\ee
and $|\det F_{\zeta}(0)|=1$. By (\ref{eqth3}) and (\ref{eqn}), we also have
$$|F_{\zeta}(\zeta)|+|F_{\overline{\zeta}}(\zeta)| \leq M^{\ast}|\det
F_{\zeta}(\zeta)|^{\frac{1}{n}} \leq\frac{M^{\ast}}{1-t|\zeta|},
$$
where $M^{\ast}=\frac{K^{\frac{1}{n}}(1+K_{1})}{K_{1}}$.

For $\zeta\in\mathbb{B}^{n}$, let $\zeta=s\theta$,
where $\theta \in \partial\IB^n$ and $s=|\zeta|$. Then
\begin{eqnarray*}
|F(\zeta)|&\leq&\int_{[0,\zeta]}|\,d F(\zeta)| =
\int_{[0,\zeta]}\big|F_{\zeta}(s\theta)\theta \,
ds+F_{\overline{\zeta}}(s\theta)\overline{\theta}\,ds\big|\\
& \leq & M^{\ast}\int_{0}^{1}\frac{ds}{1-ts} =
\frac{M^{\ast}}{t} \log \left (\frac{1}{1-t}\right ) :=M(t)  .
\end{eqnarray*}
Then by using \cite[Theorem 5]{HG1}, we have
$f(\mathbb{B}^{n})$ contains a schlicht ball with radius at least
$$B_{f}\geq\max_{0<t<1}\left\{\frac{\alpha\pi^{4n}t}{8m(4M(t))^{4n-1}}\right\},$$
where $M(t)$ is as in the statement
and $m$ is defined as in Theorem \ref{thm8}.
\qed


\end{document}